\documentclass{amsart}
\usepackage{amscd,amssymb,hyperref}
\theoremstyle{plain}
\newtheorem{prop}[subsection]{Proposition}
\newtheorem{thm}[subsection]{Theorem}
\newtheorem{cor}[subsection]{Corollary}

\theoremstyle{remark}
\newtheorem{rem}[subsection]{Remark}

\theoremstyle{definition}
\newtheorem{exm}[subsection]{Example}

\newcommand{\A}{{\mathcal A}}
\newcommand{\Ai}{{\mathcal A}_\infty}
\newcommand{\B}{{\mathcal B}}
\newcommand{\CC}{{\mathcal C}}
\newcommand{\LL}{{\mathcal L}}
\newcommand{\RR}{{\mathcal R}}

\newcommand{\Z}{{\mathbb Z}}
\newcommand{\ZN}{{\mathbb Z}_N}
\newcommand{\C}{{\mathbb C}}
\newcommand{\CP}{{\mathbb{CP}}}
\newcommand{\N}{{\mathbb N}}
\newcommand{\Zn}{{\mathbb Z}_n}
\newcommand{\Q}{{\mathbb Q}}

\newcommand{\bone}{{\mathbf 1}}
\newcommand{\bk}{{\mathbf k}}
\newcommand{\bm}{{\mathbf m}}
\newcommand{\bx}{{\mathbf x}}
\newcommand{\bl}{{\boldsymbol{\lambda}}}

\newcommand{\D}{{\Delta}}
\newcommand{\la}{{\lambda }}
\newcommand{\bul}{{\bullet }}

\renewcommand{\a}{{\alpha }}

\renewcommand{\c}{{\gamma }}
\renewcommand{\ll}{{\ell }}

\DeclareMathOperator{\rank}{rank}
\DeclareMathOperator{\ii}{i}
\DeclareMathOperator{\id}{id}
\DeclareMathOperator{\GL}{GL}

\begin{document}

\title[Cyclic covers of arrangement complements]
{Some cyclic covers of complements\\ of arrangements}
\author[D.~Cohen]{Daniel C.~Cohen$^\dag$}
\address{Department of Mathematics, Louisiana State University, 
Baton Rouge, LA 70803}
\email{\href{mailto:cohen@math.lsu.edu}{cohen@math.lsu.edu}}
\urladdr{\href{http://math.lsu.edu/~cohen/}
{http://math.lsu.edu/\~{}cohen}}
\thanks{{$^\dag$}Partially supported by 
grants LEQSF(1996-99)-RD-A-04 and LEQSF(1999-2002)-RD-A-01\\ 
from the Louisiana Board of Regents}

\author[P.~Orlik]{Peter Orlik}
\address{Department of Mathematics, University of Wisconsin, 
Madison, WI 53706}
\email{\href{mailto:orlik@math.wisc.edu}{orlik@math.wisc.edu}}

%
%
%
%
%


\begin{abstract}
Motivated by the Milnor fiber of a central arrangement, we study the 
cohomology of a family of cyclic covers of the complement of an 
arbitrary arrangement.  We give an explicit proof of the polynomial 
periodicity of the Betti numbers of the members of this family of 
cyclic covers.
\end{abstract}

\maketitle

\section{Introduction}
\label{sec:intro}

Let $\A$ be a hyperplane arrangement in $\C^\ll$, with complement 
$M(\A)=\C^\ll\setminus\bigcup_{H\in\A} H$.  The cohomology of $M(\A)$ 
with coefficients in a local system arises in a number of 
applications, both outside and inside arrangement theory.  Included 
among the former are the Aomoto-Gelfand theory of multivariable 
hypergeometric integrals \cite{AK,Gel1}, and the representation theory 
of Lie algebras and quantum groups and solutions of the 
Knizhnik-Zamolodchikov differential equation in conformal field theory 
\cite{Va}.  This note is motivated by one of the latter applications, 
the cohomology of the Milnor fiber of a central arrangement.

Let $\CC$ be a central arrangement of hyperplanes in $\C^{\ll+1}$, an 
arrangement for which each hyperplane $H\in\CC$ contains the origin.  
For each such hyperplane, let $\a_H$ be a linear form with kernel $H$.  
Then $Q=Q(\CC)=\prod_{H\in\CC}\a_H$ is a defining polynomial for the 
arrangement $\CC$, and is homogeneous of degree equal to the 
cardinality of $\CC$.  The complement $M(\CC)=\C^{\ll+1} \setminus 
Q^{-1}(0)$ may be realized as the total space of the global Milnor 
fibration
\[
F(\CC)\xrightarrow{\ \phantom{Q}\ }M(\CC)\xrightarrow{\ Q\ }\C^*,
\]
where $F(\CC)=Q^{-1}(1)$ is the Milnor fiber of $Q$, see \cite{Mi}.  
We shall refer to $F(\CC)$ as the Milnor fiber of $\CC$, and write 
$F=F(\CC)$ when the arrangement $\CC$ is understood.

Suppose that the cardinality of $\CC$ is $n$ and let 
$(x_0,x_1,\dots,x_\ll)$ be a choice of coordinates on $\C^{\ll+1}$.  
The geometric monodromy, $h:F\to F$ of the Milnor fibration is given 
by $h(x_0,\dots,x_\ll)=(\xi_{n} x_0,\dots,\xi_{n} x_\ll)$, where 
$\xi_{n}=\exp(2\pi\ii/n)$.  Since $h$ has finite order $n$, the 
algebraic monodromy $h^*:H^{q}(F;\C) \rightarrow H^{q}(F;\C)$ is 
diagonalizable and the eigenvalues of $h^*$ belong to the set of 
$n$-th roots of unity.  Denote the cohomology eigenspace of 
$\xi_{n}^k$ by $H^{q}(F;\C){_{k}}$ and write $b_q(F){_{k}} = \dim_{\C} 
H^{q}(F;\C){_{k}}$.

It is known \cite{CS1} that these cohomology eigenspaces are 
isomorphic to the cohomology of the complement of a {\em decone} of 
$\CC$ with coefficients in certain complex rank one local systems.  
See Section \ref{sec:prelim} for a summary of these results, and see
\cite{OT1} as a general reference on arrangements.  Let $\A$ be a 
decone of $\CC$, an affine arrangement in $\C^\ll$, and denote the 
rank one local systems arising in the context of the Milnor fiber by 
$\LL_k$, $1\le k \le n$.  Then $H^*(F(\CC);\C)_k \simeq 
H^*(M(\A);\LL_k)$ for each $k$.  Furthermore, the local systems 
$\LL_k$ are {\em rational} in the sense of \cite{CO}.  The results of 
this work, in the context of the Milnor fiber problem, yield 
combinatorial bounds
\begin{equation*}
\dim_\C H^q(A_\C(\A),a\wedge) \le \dim_{\C} H^{q}(F;\C){_{k}} \le 
\rank_{\Z_r} H^q(A_{\Z_r}(\A),\bar{a}\wedge), 
\end{equation*} 
where $r=n/(k,n)$ and $A_R(\A)$ is the Orlik-Solomon algebra of $\A$ 
with coefficients in the ring $R$ equipped with appropriate 
differential.  See Section \ref{sec:OS} for details.  The lower bounds 
are well known.  The upper bounds are new.

The local system $\LL_n$ is trivial, and thus corresponds to the 
constant coefficient cohomology of $M(\A)$.  This is well understood 
in terms of the Orlik-Solomon algebra.  While pursuing the remaining 
cases, we were led to a family of cyclic covers of $M(\A)$, which 
includes the Milnor fiber $F(\CC)$.  In this note, we show how a 
number of known results on the Milnor fiber extend naturally to all 
members of this family of covers, and give an explicit and elementary 
proof of the polynomial periodicity of the Betti numbers of the 
members of this family.

\section{Milnor Fibration} \label{sec:prelim}

Recall from the Introduction that $\CC$ is a central arrangement of 
$n$ hyperplanes in $\C^{\ll+1}$, with coordinates 
$(x_{0},x_{1},\dots,x_{\ll})$.  Associated with $\CC$, we have the 
defining polynomial $Q=Q(\CC)$, the complement 
$M(\CC)=\C^{\ll+1}\setminus Q^{-1}(0)$, and the Milnor fiber 
$F=F(\CC)=Q^{-1}(1)$.  The geometric monodromy $h:F\to F$ of the 
Milnor fibration has order $n=|\CC|$.  The cyclic group $\Zn$ 
generated by the geometric monodromy $h$ acts freely on $F$.  This 
free action gives rise to a regular $n$-fold covering $F \rightarrow 
F/(\Zn)$.

Consider the Hopf bundle $\C^{\ell+1}\setminus\{0\} \rightarrow 
\CP^{\ell}$ with projection map $(x_0,x_1,\dots,x_{\ell}) \mapsto 
(x_0:x_1:\dots:x_{\ell})$ and fiber $\C^*$.  Let $p:M(\CC) \rightarrow 
M^*$ denote the restriction of this projection to $M$, where $M^*= 
p(M)$.  The restriction $p_{_F} : F \rightarrow M^*$ of the Hopf 
bundle to the Milnor fiber is the orbit map of the free action of the 
geometric monodromy $h$ on $F$ and we therefore have $F/(\Zn) \cong 
M^*$.  These spaces and maps fit together with the Milnor fibration in 
the following diagram:
\[
\begin{CD}
\Zn      @>>>       F       @>{p_{_F}}>> F/\Zn
\\
@VVV                @VVV                 @VV{\cong}V \\
\C^*     @>>>       M(\CC)  @>{p}>>      M^* \\
@VVV                @VV{Q}V \\ 
\C^*/\Zn @>\cong >> \C^*  
\end{CD}
\]
Note that $M^*$ is the complement of the projective hypersurface 
defined by the homogeneous polynomial $Q$.  Thus it is the complement 
of the projective quotient of the arrangement $\CC$.  If we designate 
one of its hyperplanes, $H_\infty$, the hyperplane at infinity, the 
remaining arrangement is called the decone of $\CC$ with respect to 
$H_\infty$.  We call this $\ll$-arrangement $\A$ here and observe that 
$M(\A)=M^*$ is independent of the choice of $H_\infty$ and $|\A|=n-1$.  
We shall assume that $\A$ contains $\ell$ linearly independent 
hyperplanes.

The cohomology groups of the Milnor fiber have been studied 
extensively, see for instance \cite{OR,CS1,Ma,CScc,De}.  We summarize 
some known results from \cite{CS1}.  Since $h$ has finite order $n$, 
the algebraic monodromy $h^*:H^{q}(F;\C) \rightarrow H^{q}(F;\C)$ is 
diagonalizable and the eigenvalues of $h^*$ belong to the set of 
$n$-th roots of unity.  Denote the cohomology eigenspace of 
$\xi_{n}^k=\exp(2\pi\ii k/n)$ by $H^{q}(F;\C){_{k}}$, and denote the 
characteristic polynomial of $h^*:H^q(F;\C) \to H^q(F;\C)$ by 
$\D_q(t)=\det(t\cdot h^* - \id)$.

\begin{prop}[{\cite[1.1]{CS1}}] \label{prop:espaces}
Let $\xi_{n}^j$ and $\xi_{n}^k$ be two $n$-th roots of unity which 
generate the same cyclic subgroup of $\Zn = \langle\xi_{n}\rangle$.  
Then, for each $q$, the cohomology eigenspaces $H^{q}(F;\C){_{j}}$ and 
$H^{q}(F;\C){_{k}}$ are isomorphic.
\end{prop}

\begin{cor} \label{cor:mfcharpoly}
For each $q$, $0\le q\le \ll$, there are nonnegative integers 
$d_{k,q}$ so that
\[
\D_q(t) = \prod_{k|n} \Phi_k(t)^{d_{k,q}},
\]
where $\Phi_k(t)$ denotes the $k$-th cyclotomic polynomial.
\end{cor}

\begin{thm}[{\cite[1.6]{CS1}}] \label{thm:locsys} 
Define the rank one local system $\LL_k$ on $M(\A)$ by the 
representation $\tau_k:\pi_1(M(\A)) \to \C^*$ given by $\gamma_H 
\mapsto \xi_{n}^k$ for each meridian loop $\c_H$ about the hyperplane 
$H \in \A$.  Then, for each $k$, $1\le k \le n$, we have
\[
H^*(F(\CC);\C)_k \simeq H^*(M(\A);\LL_k).
\] 
\end{thm}

In light of this result, it is natural to study the local system on 
$M(\A)$ induced by the representation given by $\c_{H} \mapsto 
\xi_{m}^{k}$ for arbitrary $m$.  In subsequent sections, we focus on 
the context in which these local systems arise.

\section{Cyclic Covers} \label{sec:covers}

The realization of the Milnor fiber of a central arrangement as a 
cover of the complement of a decone in the previous section motivates 
the following construction.

Let $\A$ be an affine arrangement in $\C^\ll$, with coordinates 
$\bx=(x_1,\dots,x_\ll)$.  Associated with $\A$, we have the defining 
polynomial $f=Q(\A)$, and the complement $M(\A)=\C^{\ll}\setminus 
f^{-1}(0)$.  For each positive integer $m$, let $g_m:\C^* \to \C^*$ 
denote the cyclic $m$-fold covering defined by $g_m(z)=z^m$, and let 
$p_m:X_m(\A) \to M(\A)$ denote the pullback of $g_m$ along the map 
$f:M(\A)\to\C^*$, where
\[
X_m(\A) = \left\{(\bx,z) \in M(\A) \times \C^* \mid 
f(\bx)=z^m\right\}.
\]
This family of cyclic covers of $M(\A)$ generalizes the Milnor fiber 
in a number of ways which we now pursue.  First, we have the 
following.

\begin{prop} \label{prop:mfcover}
Let $\CC$ be a central arrangement of $n$ hyperplanes in $\C^{\ll+1}$, 
with defining polynomial $Q=Q(\CC)$ and Milnor fiber $F(\CC)$.  If 
$\A$ is a decone of $\CC$, then the covering spaces $p_{_F}:F(\CC) \to 
M(\A)$ and $p_n:X_n(\A) \to M(\A)$ are equivalent.
\end{prop}
\begin{proof}
The relation between the defining polynomials $f$ of $\A$ and $Q$ of 
$\CC$ is given by $Q=x_0^n\cdot f(x_1/x_0,\dots,x_n/x_0)$, and 
$F(\CC)=Q^{-1}(1)$.  Using this, it is readily checked that the map 
$X_n(\A) \to F(\CC)$ defined by $(x_1,\dots,x_\ll,z) \mapsto 
(1/z,x_1/z,\dots,x_\ll/z)$ is a homeomorphism inducing an equivalence 
of covering spaces.
\end{proof}

The characteristic homomorphism $\Phi:\pi_1(M(\A)) \to \Zn$ of the 
covering $p_{_F}:F(\CC) \to M(\A)$ was identified in \cite[1.2]{CS1}.  
It is given by $\Phi(\c_H)=g_n$, where $\c_H$ is a meridian loop about 
the hyperplane $H\in\A$ and $g_n$ is a fixed generator of $\Zn$.  A 
straightforward generalization of the proof of this fact from 
\cite[1.2]{CS1} yields

\begin{prop} \label{prop:charhom}
The characteristic homomorphism $\Phi_m:\pi_1(M(\A)) \to \Z_m$ of the 
covering $p_m:X_m(\A) \to M(\A)$ is given by $\Phi_m(\c_H)=g_m$ for a 
fixed generator $g_m$ of $\Z_m$ and meridian loops $\c_H$ about the 
hyperplanes $H$ of $\A$.
\end{prop}

The covering spaces $X_{m}(\A)$ fit together nicely in the sense of 
the following.

\begin{prop} \label{prop:subcover}
If $m=k\cdot r$, then the map $\C^\ll\times\C^* \to \C^\ll\times\C^*$ 
defined by $(\bx,z) \mapsto (\bx,z^r)$ induces a cyclic $r$-fold 
covering $p_{m,k}:X_m(\A) \to X_k(\A)$.
\end{prop}
\begin{proof}
Let $X_{k,r}(\A) \to X_k(\A)$ denote the pullback of $g_r:\C^*\to\C^*$ 
along the map $X_k(\A)\to\C^*$ defined by $(x_1,\dots,x_\ll,z) \mapsto 
z$, with
\[
X_{k,r}(\A)=\{(\bx,z,w) \in M(\A) \times \C^* \times \C^* 
\mid f(\bx)=z^k\ \text{and}\ z=w^r\}.
\]
It is then readily checked that the map $X_m(\A) \to X_{k,r}(\A)$ 
defined by $(\bx,z) \mapsto (\bx,z^r,z)$ is a homeomorphism compatible 
with the projection maps.
\end{proof}

\begin{rem} \label{rem:monodromy}
The space $X_m(\A)$ admits a self-map $h_m:X_m(\A) \to X_m(\A)$ 
defined by $h_m(\bx,z) = (\bx, \xi_m^{-1}\cdot z)$, where 
$\xi_m=\exp(2\pi\ii/m)$.  In the case $m=n$ of the Milnor fiber, the 
map $h_n:X_n(\A)\to X_n(\A)$ corresponds to the geometric monodromy 
$h:F(\CC)\to F(\CC)$ under the equivalence of covering spaces 
exhibited in the proof of Proposition~\ref{prop:mfcover}.  For 
arbitrary $m$, the ``monodromy'' map $h_m$ generates a cyclic group of 
order $m$, which acts freely on $X_m(\A)$.  The resulting regular 
$m$-fold covering $X_m(\A) \to X_m(\A)/\langle h_m\rangle$ clearly 
coincides with $p_m:X_m(\A) \to M(\A)$.  More generally, for $m=k\cdot 
r$ composite, the map $h_{m}^{k}$ generates a cyclic group of order 
$r$, which also acts freely on $X_{m}(\A)$, and the covers $X_{m}(\A) 
\to X_{m}(\A)/\langle h_{m}^{k}\rangle$ and $p_{m,k}:X_{m}(\A)\to 
X_{k}(\A)$ coincide.
\end{rem}

\section{Cohomology} \label{sec:H1}

We now study the cohomology of the covering spaces $X_{m}(\A)$.  Fix a 
basepoint $\bx_0 \in M(\A)$.  From the Leray-Serre spectral sequence 
of the fibration $p_m:X_m(\A) \to M(\A)$, we obtain $H^*(X_m(\A);\C) = 
H^*(M(\A);\LL^m)$, the cohomology of $X_m(\A)$ with trivial 
$\C$-coefficients is isomorphic to the cohomology of the base $M(\A)$ 
with coefficients in the rank $m$ local system $\LL^m$ with stalk 
$\LL^m_\bx = H^0(p_m^{-1}(\bx);\C) \simeq\C^m$.  The results presented 
in this section are natural generalizations in the context of 
arrangements of those of \cite[Section 1]{CS1}.

\begin{prop}[cf.~{\cite[1.3--1.5]{CS1}}] \label{prop:rep}
Let $T\in\GL(m,\C)$ be the cyclic permutation matrix of order $m$ 
defined by $T(\vec{e}_i)=\vec{e}_{i+1}$ for $1\le i \le n-1$ and 
$T(\vec{e}_n)=\vec{e}_1$, where $\{\vec{e}_i\}$ is the standard basis 
for $\C^m$.  Note that $T$ is diagonalizable with eigenvalues 
$\xi_m^k$, $1\le k \le m$.
\begin{enumerate}
\item The local system $\LL^m$ is induced by the representation 
$\tau^m:\pi_1(M(\A),\bx_0) \to \GL(m,\C)$ given by $\tau^m(\c_H)=T$ 
for each meridian $\c_H$.

\item The local system $\LL^m$ decomposes into a direct sum, $\LL^m = 
\bigoplus_{k=1}^m \LL^m_k$ of rank one local systems.  For each $k$, 
the local system $\LL^m_k$ is induced by the representation 
$\tau^m_k:\pi_1(M(\A),\bx_0) \to \C^*$ defined by 
$\tau^m_k(\c_H)=\xi_m^k$.

\item We have $H^*(X_m(\A);\C) = \bigoplus_{k=1}^m 
H^*(M(\A);\LL^m_k)$.
\end{enumerate}
\end{prop}

The above result provides one decomposition of the cohomology 
$H^*(X_m(\A);\C)$.  Another is given by the monodromy maps 
$h_m:X_{m}(\A)\to X_{m}(\A)$ of Remark \ref{rem:monodromy}.  Since 
$h_m$ has finite order $m$, the induced map 
$h_{m}^{*}:H^{q}(X_{m}(\A);\C) \to H^{q}(X_{m}(\A);\C)$ is 
diagonalizable, with eigenvalues among the $m$-th roots of unity.  
Denote the cohomology eigenspace of $\xi_{m}^{k}$ by 
$H^{q}(X_{m}(\A);\C)_{k}$, and let $\D^{(m)}_q(t)=\det(t\cdot h_m^* - 
\id)$ denote the characteristic polynomial of $h_m^*:H^q(X_m(\A);\C) 
\to H^q(X_m(\A);\C)$.  We then have the following generalizations of 
Theorem~\ref{thm:locsys}, Proposition~\ref{prop:espaces}, and 
Corollary~\ref{cor:mfcharpoly}.

\begin{prop}[cf.~{\cite[1.6]{CS1}}] \label{prop:decomps}
For each $k$, $1\le k \le m$, we have
\[
H^{*}(X_{m}(\A);\C)_{k} \simeq H^{*}(M(\A);\LL^m_k)
\]
\end{prop}

\begin{prop}[cf.~{\cite[1.1]{CS1}}] \label{prop:espaces2}
Let $\xi_{m}^j$ and $\xi_{m}^k$ be two $m$-th roots of unity which 
generate the same cyclic subgroup of $\Z_{m} = \langle\xi_{m}\rangle$.  
Then the cohomology eigenspaces $H^{*}(X_{m}(\A);\C){_{j}}$ and 
$H^{*}(X_{m}(\A);\C){_{k}}$ are isomorphic.
\end{prop}

\begin{cor} \label{cor:espaces3}
If $\xi_{m}^j$ and $\xi_{m}^k$ are $m$-th roots of unity which 
generate the same cyclic subgroup of $\Z_{m} = \langle\xi_{m}\rangle$, 
then $H^{*}(M(\A);\LL^m_j) \simeq H^{*}(M(\A);\LL^m_k)$.
\end{cor}

\begin{cor} \label{cor:charpoly}
For each $q$, $0\le q\le \ll$, there are nonnegative integers 
$d^{(m)}_{k,q}$ so that
\[
\D^{(m)}_q(t) = \prod_{k|n} \Phi_k(t)^{d^{(m)}_{k,q}},
\]
where $\Phi_k(t)$ denotes the $k$-th cyclotomic polynomial.
\end{cor}

We relate the cohomology of the spaces $X_m(\A)$ for various $m$ using 
these results.

\begin{thm} \label{thm:summand}
If $k$ divides $m$, then the cohomology $H^*(X_k(\A);\C)$ is a direct 
summand of $H^*(X_m(\A);\C)$.
\end{thm}
\begin{proof}
 From Proposition \ref{prop:rep}, we have 
$H^*(X_m(\A);\C)=\bigoplus_{q=1}^m H^*(M(\A);\LL^m_q)$, where the 
local system $\LL^m_q$ is induced by the representation $\tau^m_q$ 
given by $\c_H \mapsto \xi^q_m$.  Writing $m=k\cdot r$, we see that 
the representations $\tau^k_p$, $1\le p \le k$, are among the $m$ 
representations $\tau^m_q$, $1\le q \le m$.  In other words, 
$\bigoplus_{p=1}^k H^*(M(\A);\LL^k_p) = H^*(X_k(\A);\C)$ is a direct 
summand of $H^*(X_m(\A);\C)$.
\end{proof}

\begin{rem} 
This result may be interpreted in terms of the intermediate coverings 
$p_{m,k}:X_m(\A) \to X_k(\A)$ of Proposition \ref{prop:subcover} as 
follows.  One can check that the projection $p_{m,k}$ commutes with 
the monodromy maps, $p_{m,k} \circ h_m = h_k \circ p_{m,k}$.  
Consequently, the induced map $p_{m,k}^*:H^*(X_k(\A);\C) \to 
H^*(X_m(\A);\C)$ preserves the eigenspaces of $h_k^*$.  Using 
Proposition \ref{prop:decomps} and the above theorem, one can show 
that $p_{m,k}^*$ maps $H^*(X_k(\A);\C)=\bigoplus_{p=1}^k 
H^*(M(\A);\LL^k_p)$ isomorphically to the summand $\bigoplus_{p=1}^k 
H^*(M(\A);\LL^{kr}_{pr})$ of $H^*(X_m(\A);\C)$.
\end{rem}

These results also show that, to determine the cohomology of 
$X_{m}(\A)$, it suffices to compute $H^{*}(M(\A);\LL^k_{1})$ for 
divisors $k$ of $m$.  Proposition \ref{prop:decomps} and Corollary 
\ref{cor:espaces3} yield

\begin{thm} \label{thm:betti}
The Betti numbers of the space $X_{m}(\A)$ are given by 
\[
b_{q}(X_{m}(\A)) = \dim_{\C} H^{q}(X_{m}(\A);\C)= \sum_{k|m} \phi(k) 
\cdot b_q(\LL^k_1),
\]
where $\phi$ is the Euler phi function and 
$b_q(\LL^k_1)=\dim_{\C}H^{q}(M(\A);\LL^{k}_{1})$.
\end{thm}

The summand $H^{*}(M(\A);\LL^{1}_{1})$ of $H^{*}(X_{m}(\A);\C)$ 
corresponds to the constant coefficient cohomology of $M(\A)$.  This 
is well understood in terms of the Orlik-Solomon algebra defined next.

\section{Orlik-Solomon Algebra} \label{sec:OS}

Let $A=A(\A)$ be the Orlik-Solomon algebra of $\A$ generated by the 
1-dimensional classes $a_H$, $H\in \A$.  It is the quotient of the 
exterior algebra generated by these classes by a homogeneous ideal, 
hence a finite dimensional graded $\C$-algebra.  There is an 
isomorphism of graded algebras $H^*(M(\A);\C) \simeq A(\A)$.  In 
particular, $\dim A^q(\A)=b_q(\A)$ where $b_q(\A)=\dim H^q(M(\A);\C)$ 
denotes the $q$-th Betti number of $M(\A)$ with trivial local 
coefficients $\C$.  The absolute value of the Euler characteristic of 
the complement is a combinatorial invariant:
\begin{equation} \label{eq:beta}
\beta(\A)=(-1)^{\ell}\sum_{q=0}^{\ell}(-1)^q b_q(\A)=|\chi(M(\A))|.
\end{equation}
 
Let $\bl=\{\la_H\mid H \in \A\}$ be a collection of complex weights.  
Define a differential $A^q\to A^{q+1}$ by multiplication by 
$a_\bl=\sum_{H \in \A}\la_H\,a_H$.  This provides a complex 
$(A^{\bul}, a_\bl \wedge)$.  Associated to $\bl$, we have a rank one 
representation $\rho:\pi_1(M(\A))\to\C^*$ given by $\c_H\mapsto 
t_H=\exp(2\pi\ii\la_H)$ for any meridian loop $\c_H$ about the 
hyperplane $H \in\A$, and a corresponding rank one local system $\LL$ 
on $M(\A)$.  Note that $\rho$ and $\LL$ are unchanged if we replace 
the weights $\bl$ with $\bl+\bm$, where $\bm=\{m_H\mid H \in \A\}$ is 
a collection of integers.  The following inequalities are well known, 
see \cite{CO}.

\begin{prop} \label{prop:bounds}
For all $\bl$ and all $q$ we have
\[
\sup_{\bm\in\Z^{|\A|}}\dim_{\C}H^q(A^{\bul},a_{\bl+\bm}\wedge) \leq
\dim_{\C}H^q(M(\A);\LL) \leq
\dim_{\C}H^q(M(\A);\C).
\]
\end{prop}
 
For the local systems arising in the context of the covers $X_m(\A)$, 
the weights are rational.  Suppose $\la_H=k_H/N$ for all $H$, with 
integers $k_H$ and $N$, and assume without loss that the g.c.d.~of the 
$k_H$ is prime to $N$.  In this case there are better upper bounds.  
The Orlik-Solomon ideal is defined by integral linear combinations of 
the generators, hence the algebra may be defined over any commutative 
ring $R$, denoted $A_R(\A)$.  Write $A_\Q=A_\Q(\A)$.  
Left-multiplication by the element $a_\bl=\sum \la_H a_H \in A^1_\Q$ 
induces a differential on the Orlik-Solomon algebra, and we denote the 
resulting complex by $(A^{\bul}_\Q,a_\bl\wedge)$.  Similarly, 
associated to the element $a_\bk=N a_\bl=\sum k_H a_H$, we have the 
complex $(A^{\bul}_\Q,a_\bk\wedge)$.  We showed in \cite{CO} that the 
complexes $(A^{\bul}_\Q,a_\bl\wedge)$ and $(A^{\bul}_\Q,a_\bk\wedge)$ 
are chain equivalent.  The coefficients of $a_\bk$ are integers, so we 
consider the Orlik-Solomon algebra with integer coefficients and the 
associated complex $(A^{\bul}_\Z,a_\bk\wedge)$.  Let 
$(A^{\bul}_N,\bar{a}_\bk\wedge)$ be the reduction of 
$(A^{\bul}_\Z,a_\bk\wedge)\mod N$, where $A_N=A_{\ZN}$ denotes the 
Orlik-Solomon algebra with coefficients in the ring $\ZN$ and 
$\bar{a}_\bk=a_\bk\mod N$.

\begin{thm}[{\cite[4.5]{CO}}]\label{thm:modN}
Let $\bl=\bk/N$ be a system of rational weights, and let $\LL$ be the 
associated rational local system on the complement $M$ of $\A$.  Then, 
for each $q$,
\begin{equation*} \label{eq:UB}
\dim_\C H^q(M(\A);\LL) \le \rank_{\ZN} 
H^q(A^{\bul}_N,\bar{a}_\bk\wedge).
\end{equation*}
\end{thm}
There are examples in \cite{CO} which show that the inequality can be 
strict.

\section{Polynomial Periodicity} \label{sec:H2}

We continue the study of the cohomology of the spaces $X_{m}(\A)$.  We 
first investigate the implications of a well known vanishing theorem 
of Schechtman, Terao, and Varchenko in this context.  We then 
establish the polynomial periodicity of the Betti numbers of this 
family of spaces.

An edge of $\A$ is a nonempty intersection of hyperplanes and $L(\A)$ 
is the set of edges.  Given $Y\in L(\A)$, let $\A_Y=\{H \in \A\mid 
Y\subset H\}$.  Define the weight of $Y$ by $\la_Y=\sum_{H \in 
\A_Y}\la_H$.  Call $Y\in L(\A)$ {\em dense} if $\beta((\A_Y)_0)>0$ 
where $(\A_Y)_0$ is a decone of $\A_Y$.  The projective closure, 
$\Ai$, of $\A$ adds the infinite hyperplane, $H_\infty$, with weight 
$-\sum_{H \in \A}\la_H$.  Recall that $\A\subset\C^\ll$ contains $\ll$ 
linearly independent hyperplanes.

\begin{thm}[{\cite[4.3]{STV}}]
Call the local system $\LL$ nonresonant if $\la_Y \not\in \Z_{\geq 0}$ 
for every dense edge $Y\in L(\Ai)$.  In this case
\[
H^q(M(\A);\LL)=0 \mbox{~for~}q\neq \ell,\quad \text{and}\quad \dim_\C 
H^{\ell}(M(\A);\LL)=\beta(\A).
\]
\end{thm}

Recall the decomposition $H^{*}(X_{m}(\A);\C)=\bigoplus_{q=1}^{m} 
H^{*}(M(\A);\LL^{m}_{q})$ of the cohomology of $X_{m}(\A)$ from 
Proposition \ref{prop:rep}.  By Theorem \ref{thm:betti}, it suffices 
to consider the case $q=1$.  In the notation of the previous section, 
the local system $\LL^{m}_{1}$ arises from the rational and equal 
weights $\la_{H}=1/m$ for all $H\in \A$.  The projective closure, 
$\Ai$, of $\A$ is the projective quotient of the cone $\CC$ of $\A$.  
Assign the weight $-|\A|/m$ to $H_\infty$.

\begin{prop} \label{prop:vanish1}
If either {\rm (i)} $m>|\A|$, or {\rm (ii)} $|(\Ai)_Y|$ is relatively 
prime to $m$ for every dense edge $Y\in L(\Ai)$, then the local system 
$\LL^{m}_{1}$ is nonresonant.
\end{prop}
\begin{proof} 
If $Y \subset H_\infty$, then $\la_Y<0$.  Otherwise, we have 
$\la_Y=|\A_Y|/m$, which cannot be a positive integer in either case 
(i) or (ii).
\end{proof}

Call a positive integer $k$ nonresonant if $k$ satisfies either of the 
hypotheses of Proposition \ref{prop:vanish1}.  Recall the 
factorization, $\D_q^{(m)}(t)=\prod_{k|m} \Phi_k(t)^{d^{(m)}_{k,q}}$, of the 
characteristic polynomial of the monodromy $h_m^*:H^q(X_m(\A);\C) \to 
H^q(X_m(\A);\C)$ provided by Corollary~\ref{cor:charpoly}.  The 
results of Section \ref{sec:H1} and this section provide the following 
information concerning the exponents $d^{(m)}_{k,q}$ arising in this 
factorization.

\begin{prop} \label{prop:exponents}
For every $m$, we have $d^{(m)}_{1,q}=b_q(\A)$ for all $q$.  If $k$ is 
nonresonant, then $d^{(m)}_{k,q}=0$ if $q<\ll$ and 
$d^{(m)}_{k,\ell}=\beta(\A)$.
\end{prop}
\begin{proof}
The first statement follows from Proposition \ref{prop:decomps} and 
the fact that $\LL^1_1$ is the trivial local system.  The second 
statement follows from Proposition \ref{prop:vanish1}.
\end{proof}

Proposition \ref{prop:vanish1} also facilitates an elementary and 
explicit proof of the polynomial periodicity of the Betti numbers of 
the family of covering spaces $X_m(\A)$.  We refer to Sarnak and Adams 
\cite{SA,Ad} for results along these lines in greater generality, and 
to Hironaka \cite{H1,H2} and Sakuma \cite{Sk} for related results on 
branched covers of surfaces and links.  A sequence $\{a_m\}_{m \in 
\N}$ is said to be {\em polynomial periodic} if there are polynomials 
$p_1(x),\dots,p_N(x) \in \Z[x]$ so that $a_m = p_i(m)$ whenever $m 
\equiv i\mod N$.  

\begin{thm} \label{thm:pp}
For each $q$, $0\le q\le\ll$, the sequence, 
$\{b_q(X_m(\A))\}_{m\in\N}$, of Betti numbers of the cyclic covers 
$X_m(\A)$ of the complement $M(\A)$ is polynomial periodic.
\end{thm}
\begin{proof}
First note that $b_0(X_m(\A))=1$ for all $m$.

Let $N=\prod_{p\ \text{prime}} p^e$ be the product of all prime powers 
$p^e$ for which the exponent $e$ is maximal so that $p^e \le |\A|$.  
Evidently, $N$ is the smallest positive integer for which $k|N$ for 
all $k \le |\A|$.  Note also that if $m\equiv i\mod N$, then
\begin{equation} \label{eq:sets}
\{k \in \N \mid 1 \le k \le |\A|\ \text{and}\ k|m\} =  
\{k \in \N \mid 1 \le k \le |\A|\ \text{and}\ k|i\}.
\end{equation}

For $1\le q\le \ll-1$ and $1\le i \le N$, define constant 
polynomials $p_{q,i} = b_q(X_i(\A))$.  From Theorem \ref{thm:betti}, 
we have $p_{q,i} = \sum_{k|i} \phi(k)\cdot b_q(\LL^k_1)$.  Similarly, 
if $m\equiv i\mod N$, then $b_q(X_m(\A)) = \sum_{k|m} \phi(k)\cdot 
b_q(\LL^k_1)$, and by Proposition \ref{prop:vanish1}, the sum is over 
all $k \le |\A|$.  Thus, polynomial periodicity of the Betti numbers 
$b_q(X_m(\A))$ for $1\le q\le\ll-1$ follows from the relation between 
such divisors of $m$ and $i$ noted in \eqref{eq:sets} above.

The polynomial periodicity of the top Betti number $b_\ll(X_m(\A))$ 
may be established by an Euler characteristic argument as follows.  We 
have $\chi(X_m(\A))=m\cdot\chi(M(\A))$.  This, together with 
\eqref{eq:beta}, yields
\[
b_\ll(X_m(\A)) = m\cdot\beta(\A) + (-1)^{\ll+1}\left[1+ 
\sum_{q=1}^{\ll-1}(-1)^q b_q(X_m(\A))\right].
\]
Defining linear polynomials $p_{\ll,i}(x) = \beta(\A)\cdot x + 
(-1)^{\ll+1}\left[1+ \sum_{q=1}^{\ll-1}(-1)^q p_{q,i}\right]$ for each 
$i$, $1 \le i \le N$, we have $b_\ll(X_m(\A))=p_{\ll,i}(m)$ if 
$m\equiv i\mod N$.
\end{proof}

\begin{rem} 
The polynomial periodicity of the Betti numbers of more general 
classes of covers of a finite CW-complex are established in \cite{SA,Ad}.  
Noteworthy in the above proof are the explicit identifications of the 
``period'' $N$ and the polynomials $p_{q,i}(x)$ for the cyclic covers 
$X_{m}(\A)$, see the concluding remarks in \cite{H1}.  
\end{rem}

A generating function for the Betti numbers $b_q(X_m(\A))$ is given by 
the following {\em zeta function}, suggested by A.~Adem.  For each 
$q$, $0 \le q \le \ll$, define
\begin{equation} \label{eq:zeta}
\zeta_{\A,q}(s) = \sum_{m=1}^{\infty} \frac{b_{q}(X_{m}(\A))}{m^{s}}.
\end{equation}

\begin{thm} \label{thm:zeta}
We have
\begin{equation*}
\zeta_{\A,q}(s) =
\zeta(s) \cdot \left[\sum_{k \le |\A|}\frac{\phi(k)\cdot 
b_q(\LL^k_1)}{k^s}
+ \delta_{q,\ll} \cdot \beta(\A) 
\sum_{k > |\A|}\frac{\phi(k)}{k^s}\right]
\end{equation*}
where $\zeta(s)$ is the classical Riemann zeta function and 
$\delta_{q,\ll}$ is the Kronecker delta.
\end{thm}
\begin{proof}
 From Theorem \ref{thm:betti}, we have $b_q(X_m(\A))=\sum_{k|m} 
\phi(k)\cdot b_q(\LL^k_1)$.  If $k>|\A|$, we have $b_q(\LL^k_1)=0$ 
for $q<\ll$ and $b_\ll(\LL^k_1)=\beta(\A)$ by Proposition 
\ref{prop:vanish1}.  A calculation using these observations yields 
the result.
\end{proof}

\section{Bounds and Examples}
\label{sec:examples}

The results of Sections \ref{sec:H1} and \ref{sec:H2} show that, to 
determine the cohomology of $X_m(\A)$, it suffices to compute 
$H^*(M(\A);\LL^k_1)$ for those divisors $k$ of $m$ for which $k<|\A|$ 
and hypothesis (ii) of Proposition \ref{prop:vanish1} fails.  
Combining the results of Proposition \ref{prop:bounds} and Theorem 
\ref{thm:modN}, we have combinatorial bounds on the local system Betti 
numbers,
\begin{equation} \label{eq:bounds}
\sup_{\bm\in\Z^{|\A|}}\dim_{\C}H^q(A^{\bul},a_{\bl+\bm}\wedge) \le 
b_q(\LL^k_1) \le \rank_{\Z_k} 
H^q(A^{\bul}_{\Z_k},\bar{a}_\bone\wedge),
\end{equation}
where $k\cdot \bl = \bone$ and $a_\bone=\sum_{H\in\A} a_H$.  
Evidently, if the two extreme non-negative integers in the above 
inequalities are equal, then $b_q(\LL^k_1)$ is determined.  In 
particular, if $\rank_{\Z_k} H^q(A^{\bul}_{\Z_k}, 
\bar{a}_\bone\wedge)=0$, we have $b_q(\LL^k_1)=0$ as well.  We 
conclude with several examples which illustrate the utility of these 
bounds.

\begin{exm} \label{exm:maclane}
Let $\CC$ be a realization of the MacLane ($8_{3}$) configuration, 
with defining polynomial 
$Q(\CC)=xy(y-x)z(z-x-\xi_3^{2}y)(z+\xi_3y)(z-x)(z+\xi_3^{2}x+\xi_3y)$, 
and let $\A$ be a decone of $\CC$.  The Poincar\'e polynomial of $\A$ 
is $P(\A,t)=1+7t+13t^2$, and $\beta(\A)=7$.  A calculation in the 
Orlik-Solomon algebra of $\A$ reveals that $H^q(A^{\bul}_{\Z_k}, 
\bar{a}_\bone\wedge)=0$ for $q\neq 2$ and all $k>1$.  Thus, 
$P(X_{m}(\A),t)=1+7t+(6+7m)t^2$ for all $m$.  In particular, for 
$m=8=|\CC|$, the Poincar\'e polynomial of the Milnor fiber of the 
MacLane arrangement is $P(F(\CC),t)=1+7t+62t^2$.
\end{exm}

\begin{exm} \label{exm:selberg}
Let $\A$ be the Selberg arrangement in $\C^{2}$, with defining 
polynomial $Q(\A)=xy(x-y)(x-1)(y-1)$.  The Poincar\'e polynomial of 
$\A$ is given by $P(\A,t)=\sum_{q\ge 0} b_q(\A)t^q = 1+5t+6t^2$, and 
we have $\beta(\A)=2$, see \eqref{eq:beta}.  The dense edges of $\Ai$ 
all have cardinality $3$, so by Proposition \ref{prop:vanish1}, if 
$k>5$ or $k$ is prime to $3$, the local system $\LL^k_1$ on $M(\A)$ is 
nonresonant.  Thus if $(m,3)=1$, the Poincar\'e polynomial of the 
cover $X_m(\A)$ is $P(X_{m}(\A),t)=1+5t+(4+2m)t^2$.

If $k=3$, one can check that $\dim_{\C}H^q(A^{\bul},a_{\bl+\bm}\wedge) 
= \rank_{\Z_3} H^q(A^{\bul}_{\Z_3},\bar{a}_\bone\wedge) = 1$, where 
$3\bl=\bone$ and $\bm=(\dots m_{H} \dots)\in\Z^{5}$ satisfies 
$m_{H}=-1$ if $H=\{x-y=0\}$ and $m_{H}=0$ otherwise.  Consequently, 
$b_{1}(\LL^{3}_{1})=1$ as well, and if $3$ divides $m$, we have 
$P(X_{m}(\A),t)=1+7t+(6+2m)t^2$.  It follows that the zeta function 
$\zeta_{\A,1}(s)$ of \eqref{eq:zeta} is given by 
$\zeta_{\A,1}(s)=\zeta(s)\cdot [5+2\cdot 3^{-s}]$.

Since the Selberg arrangement is a decone of the braid arrangement 
$\B$ of rank three, the cover $X_{6}(\A)$ is homeomorphic to the 
Milnor fiber $F(\B)$, and $P(F(\B),t)=1+7t+18t^{2}$, as is well known.  
For further calculations along these lines, see \cite{CScc,De}.
\end{exm} 

\begin{exm}\label{exm:hessian} 
Let $\CC$ be the Hessian configuration, with defining polynomial 
$Q(\CC)=x_1x_2x_3\prod_{i,j=0,1,2} (x_1+\xi_{3}^i x_2+\xi_{3}^j x_3)$, 
and let $\A$ be a decone of $\CC$ with Orlik-Solomon algebra $A$.  The 
Poincar\'e polynomial of $\A$ is $P(\A,t)=1+11t+28t^{2}$, and 
$\beta(\A)=18$.  One can check that $H^q(A^{\bul}_{\Z_k}, 
\bar{a}_\bone\wedge)=0$ for $q\neq 2$ if $k \neq 2,4$, and that, if 
$k=2,4$,
\[
\rank_{\Z_k} H^q(A^{\bul}_{\Z_k},\bar{a}_\bone\wedge)=
\begin{cases}
2&\text{if $q=1$,}\\
20&\text{if $q=2$,}\\
0&\text{otherwise.}
\end{cases}
\]
So $b_{1}(\LL^{k}_{1})=0$ and $b_{2}(\LL^{k}_{1})=18$ if $k\neq 2,4$, 
while $b_{1}(\LL^{k}_{1})\le 2$ and $b_{2}(\LL^{k}_{1})\le 20$ if $k= 
2,4$.

Concerning the lower bound of \eqref{eq:bounds}, it is known that the 
resonance variety $\RR_{1}(\CC)$ of the Hessian arrangement has a 
non-local three-dimensional component $S\subset\C^{12}$, see 
\cite[5.8]{CS4}, \cite[3.3]{Li}.  For $\bl \in S$, we have $\dim_\C 
H^1(A^{\bul}(\CC),a_\bl\wedge)=2$, see \cite[3.12]{Fa}.  For an 
appropriate ordering of the hyperplanes of $\CC$, this component has 
basis
\[
\vec{e}_5+\vec{e}_7+\vec{e}_{12}-\vec{e}_1-\vec{e}_2-\vec{e}_3,
\hfill 
\vec{e}_4+\vec{e}_9+\vec{e}_{11}-\vec{e}_1-\vec{e}_2-\vec{e}_3,
\hfill 
\vec{e}_6+\vec{e}_8+\vec{e}_{10}-\vec{e}_1-\vec{e}_2-\vec{e}_3. 
\]
Using this basis, one can show that, for $a_\bl = \frac{1}{k}a_\bone 
\in \C^{12}$, there exists $\bm\in\Z^{12}$ so that ${\bl+\bm}\in S$ if 
and only if $k=2,4$.  Thus, 
$\sup_{\bm\in\Z^{12}}\dim_{\C}H^1(A^{\bul}(\CC),a_{\bl+\bm}\wedge)=2$ 
for these $k$.  A standard argument then shows that 
$\sup_{\bm\in\Z^{11}}\dim_{\C}H^1(A^{\bul}(\A),a_{\bl+\bm}\wedge)=2$ 
for $a_\bl = \frac{1}{k}a_\bone \in \C^{11}$ as well.  Consequently, 
the inequalities in the upper bounds on $b_{q}(\LL^{k}_{1})$ noted 
above are in fact equalities.  It follows that the zeta function 
$\zeta_{\A,1}(s)$ is given by 
$\zeta_{\A,1}(s)=\zeta(s)\cdot[11+2\cdot 2^{-s}+4\cdot4^{-s}]$.

These calculations determine the Betti numbers of the cover $X_m(\A)$ 
for any $m$, as well as the dimensions of the eigenspaces of the maps 
$h_m^*:H^*(X_m(\A);\C)\to H^*(X_m(\A);\C)$.  In particular, for 
$m=12=|\CC|$, the Poincar\'e polynomial of the Milnor fiber of the 
Hessian arrangement is $P(F(\CC),t)=1+17t+232t^2$, and the 
characteristic polynomials $\D_q(t)=\D^{(12)}_{q}$ of the algebraic 
monodromy are
\[
\D_0(t)=t-1,\  \D_1(t)=(t-1)^{9}(t^4-1)^{2},\ \text{and}\  
\D_2(t)=(t-1)^{8}(t^4-1)^2(t^{12}-1)^{18}.
\]
\end{exm}

\begin{exm}\label{exm:ceva}
Let $\A$ be the Ceva(3) arrangement in $\C^3$, with defining 
polynomial $Q(\A)=(x^3-y^3)(x^3-z^3)(y^3-z^3)$.  It is known that 
$\dim_\C H^1(A^{\bul},a_{\bl}\wedge) \le 1$ for all $\bl$.  In the 
case $k=3$, is is also known that $b_1(\LL^3_1)=\rank_{\Z_3} 
H^1(A^{\bul}_{\Z_3},\bar{a}_\bone\wedge)=2$.  This example is 
discussed in detail in \cite[4.5]{Fa}, \cite[3.3]{Li}, 
\cite[6.2]{CS4}, and \cite[3.5 and 4.7]{CO}.

As illustrated by this arrangement, the lower bound of 
\eqref{eq:bounds} may be strict.  On the other hand, we know of no 
example where the upper bound of \eqref{eq:bounds} is strict.
\end{exm}

\bibliographystyle{amsalpha}

\begin{thebibliography}{STV}

\bibitem[Ad]{Ad} S. Adams,
{\em Representation varieties of arithmetic groups and 
polynomial periodicity of Betti numbers}, 
Israel J. Math. 88 (1994), 73--124. 

\bibitem[AK]{AK} K.~Aomoto, M.~Kita, {\em Hypergeometric Functions} 
(in Japanese), Springer-Verlag, 
1994.

\bibitem[CO]{CO} D. Cohen, P. Orlik, 
{\em Arrangements and local systems}, preprint, 1999;
\texttt{\href{http://xxx.lanl.gov/abs/math.AG/9907117}
{math.AG/9907117}}.

\bibitem[CS1]{CS1} D.~Cohen, A.~Suciu, {\em On Milnor fibrations of 
arrangements}, J. London Math.  Soc.  \textbf{51} (1995), 105--119.

\bibitem[CS2]{CScc} \bysame,
{\em Homology of iterated 
semidirect products of free groups}, J. Pure Appl.  Algebra 
\textbf{126} (1998), 87--120.

\bibitem[CS3]{CS4} \bysame, {\em Characteristic varieties of 
arrangements}, Math. Proc. Cambridge Philos. Soc. \textbf{127} 
(1999), 33--54.

\bibitem[De]{De} G.~Denham, {\em Local systems on the complexification 
of an oriented matroid}, Thesis, University of Michigan, 1999.

\bibitem[Fa]{Fa} M. Falk, {\em Arrangements and cohomology}, Ann.  
Comb. \textbf{1} (1997), 135--157.

\bibitem[Ge]{Gel1} I. M. Gelfand, 
{\em General theory of hypergeometric functions},
Soviet Math. Dokl. \textbf{33} (1986), 573-577.

\bibitem[H1]{H1} E.~Hironaka,
{\em Polynomial periodicity for Betti numbers of covering surfaces},
Invent. Math. \textbf{108} (1992), 289--321.

\bibitem[H2]{H2} \bysame,
{\em Intersection theory on branched covering surfaces 
and polynomial periodicity},
Internat. Math. Res. Notices (1993), 185--196.  

\bibitem[Li]{Li} A.~Libgober, 
{\em Characteristic varieties of algebraic curves}, preprint, 1998;
\texttt{\href{http://xxx.lanl.gov/abs/math.AG/9801070}
{math.AG/9801070}}.

\bibitem[Ma]{Ma} D.~Massey,  
{\em Perversity, duality and arrangements in $\C^3$}, 
Topology Appl. \textbf{73} (1996), 169--179. 

\bibitem[Mi]{Mi} J.~Milnor, {\em Singular Points of Complex 
Hypersurfaces}, Annals of Math. Studies \textbf{61}, Princeton 
University Press, 1968.

\bibitem[OR]{OR} P.~Orlik, R.~Randell,
{\em The Milnor fiber of a generic arrangement}, 
Arkiv f\"{u}r Mat. \textbf{31} (1993), 71--81.

\bibitem[OT]{OT1} P.~Orlik, H.~Terao, {\em Arrangements of 
Hyperplanes}, Grundlehren Math.  Wiss., vol.~300, Springer-Verlag, 
1992.

\bibitem[Sk]{Sk} M. Sakuma,
{\em Homology of abelian coverings of links and spatial graphs},
Canad. J. Math. \textbf{47} (1995), 201--224.

\bibitem[SA]{SA} P. Sarnak, S. Adams,
{\em Betti numbers of congruence groups},
Israel J. Math. \textbf{88} (1994), 31--72.

\bibitem[STV]{STV} V.~Schechtman, H.~Terao, A.~Varchenko, {\em 
Cohomology of local systems and the Kac-Kazhdan condition for 
singular 
vectors}, J. Pure Appl.  Algebra \textbf{100} (1995), 93--102.

\bibitem[Va]{Va} A.~Varchenko, {\em Multidimensional Hypergeometric 
Functions and Representation Theory of Lie Algebras and Quantum 
Groups}, Adv. Ser. Math. Phys., vol.  21, World 
Scientific, 
1995.

\end{thebibliography}

\end{document}